\documentclass[12pt]{article}
\usepackage{amsmath}
\usepackage{amssymb}
\usepackage{amsthm}
\usepackage{cite}
 \usepackage{longtable}
\usepackage[latin1]{inputenc}

\newtheorem{definition}{Definition}

\newtheorem{theorem}{Theorem}

\def\PG{\mathrm{PG}}

\begin{document}

\title{Classification of minimal 1-saturating sets in $\PG(v, 2), 2 \leq v \leq 6$}
\date{}
\maketitle

{\large\emph{Alexander A. Davydov$\,^1$, Stefano Marcugini$\,^2$, Fernanda Pambianco$\,^2$}}\\
{\footnotesize $^1$Kharkevich Institute for Information Transmission Problems,
Russian Academy of  Sciences, Bol'shoi Karetnyi pereulok\ 19, Moscow, 127051, Russian Federation} \\
{\footnotesize$^2$Dipartimento di Matematica e Informatica,
Universit\`{a} degli Studi di Perugia, Via Vanvitelli~1, Perugia, 06123, Italy}\\
{\footnotesize E-mails: adav@iitp.ru, stefano.marcugini@unipg.it, fernanda.pambianco@unipg.it} \\
\medskip

\emph{
\textbf{Abstract:} The classification of all the minimal $1$-saturating sets in $%
\PG(v,2)$ for $2 \leq v \leq 5$, and the classification of the smallest and of the second smallest minimal $1$-saturating sets in $\PG(6,2)$ are presented. These results have been found using a computer-based exhaustive search. \\}
\medskip

\emph{
}

\section{Introduction}

Let $F_{q}$ be the Galois field of $q$ elements and let $\PG(v,q)$ be the $v$
-dimensional projective space over $F_{q}$. For an introduction to
geometrical objects in such spaces, see~\cite{Hirs, HirsSt}.

For an integer $\varrho $ with $0\leq \varrho \leq n$ we say that a set of
points $S\subseteq PG(v,q)$ is \linebreak $\varrho ${\em -saturating} if for
any point $x\in PG(v,q)$ there exist $\varrho +1$ points in $S$ generating a
subspace of $PG(v,q)$ in which $x$ lies and $\varrho $ is the smallest value
with such property, cf. \cite{DMP, DavO2, Ughi}.

Note that the term ``saturated'' for points in $S$ was applied in \cite{Ughi}
and then was used in some papers. But in \cite{PamSt} the points of $%
PG(v,q)\setminus S$ are said to be saturated and this seems to be more
natural. Therefore in \cite{DMP, DavO2} and here the points in $S$ are
called ``saturating''. 
In \cite{G2007}, see also the references therein, saturating sets
are called "dense sets".
Note also that in \cite{BPW} saturating sets are
called ``$R$-spanning sets''. Finally, in some works the points of $%
PG(v,q)\setminus S$ are called to be ``covered''. This term seems acceptable
too.

A $\varrho $-saturating set of $k$ points is called {\em minimal} if it does
not contain a $\varrho $-saturating set of $k-1$ points \cite{DMP, Ughi}.

In this paper we consider minimal 1-saturating sets in binary projective
spaces $PG(v,2)$. A set $S\subset PG(v,2)$ is 1-saturating if any point of $%
PG(v,2)\setminus S$ lies on a bisecant of $S$.

Arcs in $PG(2,2)$ and caps in $PG(v,2),$ $v\geq 3,$ are sets of points, no
three of which are collinear. Complete arcs and caps are minimal
1-saturating sets \cite{DMP, Ughi} which we call ``CA sets'' for
complete arcs and ``CC sets'' for complete caps. For sizes, constructions,
and estimates of binary CA and CC sets, see, for example, \cite{ BW1, BW2, DavPrep02, DFP, DGMP2010, DMPbin, DMPJG, DT, FPJG, GDT, GL2010, Hirs, HirsSt, KhaLis, Lison, Weh1},
and the references therein.

On the other hand, a minimal 1-saturating set may contain three points of
the same line. Then it is neither an arc nor a cap. We call such minimal
1-saturating set an ``NA set'' in $PG(2,2)$ and an ``NC set'' in $PG(v,2),$ $%
v\geq 3$. NC sets have a more wide spectrum of possible sizes than CC sets.
Some constructions, sizes, and estimates for binary NC sets are given in 
\cite{BDGMP2016, Handbook, BPW, Coh, DMP, FG2003, GDT, G2007, G2013, GL2010, GrSl, KaikRos, Ughi}, and the references therein, either
directly or they can be obtained from those for $q = 2$.
Of particular interest is \cite{DMPbin}, where several constructions of minimal 1-saturating sets in 
binary projective spaces $PG(v,2)$ are presented.
In \cite{GL2010}, the authors observe  that a minimal 1-saturating set can be obtained from a complete cap $S$ by fixing some $s \in S$ and replacing every point $s' \in S \setminus  \{s\}$ by the third point on the line through $s$ and $s'$.  From here on, we will denote this construction by GL.

For NC sets we can use results of the linear covering codes theory, e.g., of
\cite{Handbook, BPW, Coh, DMPbin, DFMP-IEEE-LO, GDT, GrSl, KaikRos}%
, due to the following considerations. A $q$-ary linear code with
codimension $r$ has {\em covering radius} $R$ if every $r$-positional $q$%
-ary column is equal to a linear combination of $R$ columns of a parity
check matrix of this code and $R$ is the smallest value with such property.
For an introduction to coverings of vector spaces over finite fields and to
the concept of code covering radius, see \cite{Handbook, Coh}. The
points of a $\varrho $-saturating $n$-set in $PG(r-1,q)$ can be considered
as columns of a parity check matrix of a $q$-ary linear code of length $n,$
codimension $r,$ and covering radius $\varrho +1.$ This correspondence is
remarked and used in many works, see, for example, \cite{BPW, DMP, DavO2}, and the references therein.

For given codimension and covering radius, the linear covering codes theory
\cite{Handbook, Coh}, is interested in codes of the smallest length since
they have small covering density. In a geometric perspective, saturating
sets of the smallest size are also interesting as extremal objects.


In terms of linear covering codes, the concept of minimal saturating sets corresponds to the concept of locally optimal linear covering code; see \cite{DFMP-IEEE-LO}. A locally optimal code is nonshortening in the sense that one cannot remove any column from a parity-check matrix without increasing the code covering radius.

At present minimal saturating sets seem to be studied insufficiently. In
general, their smallest sizes and the spectrum of possible sizes are
unknown. Relatively a few constructions of minimal saturating sets are
described in literature.

Note that in $\PG(v,2),$ a complete cap of maximal size is the complement of a hyperplane, see \cite{DT}, and its stabilizer group is  $  ASL(v,2) $, while a minimal 1-saturating set of maximal size that is not a cap is a hyperplane together with a point outside it, see   \cite[Corollary 1]{DMP}, and  its stabilizer group is  $ PSL(v,2) $. In  both the cases the size of the set is $2^v$. \\
The Structure Theorem of Davydov and Tombak gives a characterization of "large" binary caps:
\begin{theorem}[\cite{DT}]
Any "large" (cardinality $\geq 2^{v-1}+2$) complete cap in $\PG(v,2)$ is obtained by a repeated application of the doubling construction to a "critical" complete cap (cardinality $2^{k-1}+1$) in $\PG(k,2)$ for some $k < v$. 
\end{theorem}

In \cite{GL2010} it is stated that in $\PG(v,2)$ every 1-saturating set of size at least $\frac{11}{36} 2^{v+1} +3$ either is a complete cap or can be obtained from a complete cap $S$ by construction GL, that the 1-saturating sets of the second largest size are the complete cap of size $5 \times 2^{n-3}$ and the corresponding NC set defined as above, and that the  third largest size is smaller than $\frac{11}{36} 2^{v+1} +3$. 

Note that by applying construction GL to the complement of a hyperplane, you obtain a hyperplane and a point ouside it, while by applying it to the complete cap of size five in $\PG(3,2)$, you obtain a projectively equivalent complete cap; the same happens by applying construction GL to the complete cap of size 17 in $\PG(5,2)$ whose stabilizer group has order 40320.

In this paper we present the classification of all the minimal $1$-saturating sets in $%
\PG(v,2)$ for $2 \leq v \leq 5$, and the classification of the smallest and of the second smallest minimal $1$-saturating sets in $\PG(6,2)$, giving for each set the list of its points, the description of its stabilizer group, and a reference to a theoretical construction when it is known. This classification has been obtained by computer.

A summary of these results appeared for the first time in \cite[Section 5]{DMPbin},  where the structure of a minimal 1-saturating 19-set in $\PG(6,2)$ is also described in detail.  

\section{Classification of minimal 1-saturating sets in $\PG(v,2)$, $2 \leq v \leq 5$ and of small minimal 1-saturating sets in $\PG(6,2)$ }\label{classification_sec}

We obtained the classification of the minimal 1-saturating sets in $\PG(v,2), 2 \leq v \leq 5$ and of the small minimal 1-saturating sets in $\PG(6,2)$  using an exhaustive computer search based on a backtracking algorithm \cite{DMP}.
The algorithm exploits equivalence properties among sets of points of $\PG(v,2),$ to reduce the search space.
However several projectively equivalent copies of the same minimal 1-saturating set can be obtained. Therefore the examples have been classified using MAGMA; see \cite{magma}.
 Using Magma,  the stabilizer group has been computed and identified, if not too big. Then the names of the groups have been determined using GAP; see \cite{GAP}. The structure of the stabilizer group of the complete caps obtained by \cite[Construction D]{DMPbin} is described in \cite{DMP-Doubl2017OC}.

In Table 1 we give the summary of the complete classification of minimal 1-saturating $k$%
-sets in $\PG(v,2),$ $v\leq 5,$ for all $k$, and in $\PG(6,2)$ for $k\leq 20$.
For ``type'' CA, CC, NA, and NC, see Introduction. The notation $n$ means
the number of objects of type noted. ``Stab. group'' gives either the order
of the stabilizer group if $n=1$ or the interval of the orders if $n>1.$
 Table 1 appeared for the first time in \cite[Section 5]{DMPbin}.

\begin{definition} $\;$\\
Let $t_{2}(v,q)$ be the smallest size of a complete arc in $\PG(2,q)$ and the smallest size of a complete cap in $\PG(v,q), v \geq 3$.\\
Let $\ell(v,q,1)$ be the smallest size of a minimal 1-saturating set  in $\PG(v,q)$.\\
Let $m(v,q,1)$ be the greatest size of a minimal 1-saturating set  in $\PG(v,q)$.\\
Let $m^{\prime}(v,q,1)$ be the second greatest size of a minimal 1-saturating set  in $\PG(v,q)$.\\
Let $m^{\prime \prime}(v,q,1)$ be the third greatest size of a minimal 1-saturating set  in $\PG(v,q)$.
\end{definition}

By Table 1, we have
\begin{eqnarray}
t_{2}(2,2) &=&\ell(2,2,1)=m^{\prime \prime }(2,2,1)=m^{\prime
}(2,2,1)=m(2,2,1)=4.  \nonumber \\
t_{2}(3,2) &=&\ell(3,2,1)=m^{\prime \prime }(3,2,1)=5,\text{ }m^{\prime
}(3,2,1)=6,\text{ }m(3,2,1)=8.  \nonumber \\
t_{2}(4,2) &=&\ell(4,2,1)=9,\text{ }m^{\prime \prime }(4,2,1)=10,\text{ }%
m^{\prime }(4,2,1)=11,\text{ }m(4,2,1)=16.  \nonumber \\
t_{2}(5,2) &=&\ell(5,2,1)=13,\text{ }m^{\prime \prime }(5,2,1)=18,\text{ }%
m^{\prime }(5,2,1)=20,\text{ }m(5,2,1)=32.  \nonumber \\
\ell(6,2,1) &=&19.\quad t_{2}(6,2)=21.  \label{computres}
\end{eqnarray}

  \begin{footnotesize}
   \begin{longtable}
{|c|c|c|c|c||c|c|c|c|c|}

  \caption{ Complete classification of minimal
1-saturating $k$-sets in $PG(v,2)$, $v\leq 5$, for all $k$, and
in $PG(6,2)$ for $k\leq 20$}

\endfirsthead
\multicolumn{8}{r}{\textit{(The table continues in the next page)}}
\endfoot
\multicolumn{5}{l}{Table 1  continue}
\endhead

\endlastfoot

\hline
\rule[0.1 mm]{0mm}{5 mm}
v & k & \text{Type} & n & \text{Stab. group} & v & k & \text{Type} & n &
\text{Stab. group} \\ \hline
2 & 4 &
\rule[1.5 mm]{0mm}{5 mm}
$
\begin{array}{c}
\text{CA} \\
\text{NA}
\end{array}
$
&
$
\begin{array}{c}
1 \\
1
\end{array}
$
&
$
\begin{array}{c}
24 \\
6
\end{array}
$
& 5 & 14 & \text{NC} & 19 & 8\ldots 56448 \\ \hline
\rule[.2 mm]{0mm}{5 mm}
3 & 5 & \text{CC} & 1 & 120 & 5 & 15 & \text{NC} & 14 & 4\ldots 72 \\ \hline
\rule[.2 mm]{0mm}{5 mm}
3 & 6 & \text{NC} & 1 & 72 & 5 & 16 & \text{NC} & 15 & 2\ldots 12 \\ \hline
\rule[1.5 mm]{0mm}{5 mm}
3 & 8 &
$
\begin{array}{c}
\text{CC} \\
\text{NC}
\end{array}
$
&
$
\begin{array}{c}
1 \\
1
\end{array}
$
&
$
\begin{array}{c}
1344 \\
168
\end{array}
$
& 5 & 17 &
$
\begin{array}{c}
\text{CC} \\
\text{NC}
\end{array}
$
&
$
\begin{array}{c}
5 \\
48
\end{array}
$
&
$
\begin{array}{c}
384\ldots 40320 \\
2\ldots 8064
\end{array}
$

\\ \hline
4 & 9 &
\rule[1.5 mm]{0mm}{5 mm}
$
\begin{array}{c}
\text{CC} \\
\text{NC}
\end{array}
$
&
$
\begin{array}{c}
1 \\
1
\end{array}
$
&
$
\begin{array}{c}
336 \\
144
\end{array}
$
& 5 & 18 &
$
\begin{array}{c}
\text{CC} \\
\text{NC}
\end{array}
$
&
$
\begin{array}{c}
1 \\
108
\end{array}
$
&
$
\begin{array}{c}
10752 \\
2\ldots 120960
\end{array}
$
\\ \hline
4 & 10 &
\rule[1.5 mm]{0mm}{5 mm}
$
\begin{array}{c}
\text{CC} \\
\text{NC}
\end{array}
$
&
$
\begin{array}{c}
1 \\
6
\end{array}
$
&
$
\begin{array}{c}
1920 \\
8\ldots 1008
\end{array}
$
& 5 & 20 &
$
\begin{array}{c}
\text{CC} \\
\text{NC}
\end{array}
$
&
$
\begin{array}{c}
1 \\
1
\end{array}
$
&
$
\begin{array}{c}
184320 \\
9216
\end{array}
$
\\ \hline
\rule[1.5 mm]{0mm}{5 mm}
4 & 11 & NC & 1 & 10 & 5 & 32 &
$
\begin{array}{c}
\text{CC} \\
\text{NC}
\end{array}
$
&
$
\begin{array}{c}
1 \\
1
\end{array}
$
&
$
\begin{array}{c}
\end{array}
$
\\ \hline
4 & 16 &
\rule[1.5 mm]{0mm}{5 mm}
$
\begin{array}{c}
\text{CC} \\
\text{NC}
\end{array}
$
&
$
\begin{array}{c}
1 \\
1
\end{array}
$
&
$
\begin{array}{c}
322560 \\
20160
\end{array}
$
& 6 & 19 & \text{NC} & 5 & 32\ldots 5760 \\ \hline
5 & 13 &
\rule[1.5 mm]{0mm}{5 mm}
$
\begin{array}{c}
\text{CC} \\
\text{NC}
\end{array}
$
&
$\begin{array}{c}
1 \\
7
\end{array}
$
&
$
\begin{array}{c}
1152 \\
32\ldots 4032
\end{array}
$
& 6 & 20 & \text{NC} & 36 & 4\ldots 2880 \\ \hline
\end{longtable}
   \end{footnotesize}

\smallskip

The relation $t_{2}(6,2)=21$ is based on the facts that in $PG(6,2)$ there
is a complete 21-cap \cite[Th. 3]{GDT} but there are not complete $k$-caps
with $k\leq 20,$ see Table 1 and \cite{KhaLis}. Note also that in
\cite[p. 222]{GDT} the conjecture was done that this relation holds.

The values of $\ell(v,2,1),$ $v\leq 6,$ and $t_{2}(v,2),$ $v\leq 5,$ are given
also in \cite[Table 2]{Handbook} and\linebreak\ \cite[Tables 3.1,4.2]{FPJG},
respectively.
The classification of the complete caps in $PG(v,2), v \leq 6$ can be found  in \cite{KhaLis}; in  \cite{BMP} the classification of all caps, complete and incomplete in PG(5,2) is given, together with the list of the points and the description of the stabilizer group.

In \cite[Remark 5, p. 271]{DT} five distinct complete 17-caps in $PG(5,2)$
are constructed and the conjecture is done that other nonequivalent 17-caps
in $PG(5,2)$ do not exist. This conjecture is proved by an exhaustive
computer search in  \cite{DMPbin} (see Table 1, $k=17,$ type~CC) and in~\cite
{KhaLis}. This fact allows us to obtain all nonequivalent
complete $17\cdot 2^{v-5}$-caps in $PG(v,2),$ $v\geq 6,$ by $(v-5)$-fold
applying Construction DC to a complete 17-cap in $PG(5,2)$ \cite{DT}. 
Note that the five complete 17-caps in $PG(5,2)$ can by obtained using the construction described in \cite[Theorem 2.4]{KhaLis}; two of them can be obtained also using the construction $L_{21}$ of \cite{DMPbin}. 

\newpage

The following tables give the classification of all the minimal 1-saturating sets in $\PG(v,2)$, $2 \leq v \leq 5$ and of the smallest and the second smallest minimal 1-saturating sets in $\PG(6,2)$. We denote a point $P$ of $\PG(v,2)$ by the decimal integer of which $P$ is the binary representation. When the order $i$  of a stabilizer group is too big to be identified by Magma, we denoted the group as $G_i$. When possible, we indicate the construction giving the example: KL denotes \cite[Theorem 2.4]{KhaLis}, see  \cite{DMPbin} for the other symbols.

   \begin{footnotesize}
\begin{longtable}
 {|c|c|c|c|c|}
 \caption{Classification of the minimal 1-saturating sets in $\PG(2,2)$}

\endfirsthead
\multicolumn{5}{r}{\textit{(The table continues in the next page)}}
\endfoot
\multicolumn{5}{l}{ Table 2  continue}
\endhead

\endlastfoot

   \hline
\rule[-1 mm]{0mm}{5 mm}
    Size  & Type   &
  \begin{minipage}{2.3 cm}
 \smallskip
  \centering
  List of \\
 points \\
 \medskip
  \end{minipage}   &
  \begin{minipage}{2.3 cm}
 \smallskip
  \centering
  Stabilizer\\
 group \\
 \medskip
  \end{minipage}   &
   \begin{minipage}{0.9 cm}
 \smallskip
  \centering
 Cons-\\
truc-\\
 tion \\
 \medskip
  \end{minipage}   \\
   \hline
\rule[-1 mm]{0mm}{5 mm}
$ 4 $ & CA & $ \{1, 2, 4, 7\}  $ & $ ASL(2,2) $ & H  \\
\hline
\rule[-1 mm]{0mm}{5 mm}
$ 4 $ & NA & $ \{1, 2, 4, 6\}  $ & $ PSL(2,2) $ &  A \\
\hline
  \end{longtable}
   \end{footnotesize}


 \begin{footnotesize}
   \begin{longtable}
 {|c|c|c|c|c|}

  \caption{Classification of the minimal 1-saturating sets in $\PG(3,2)$}

\endfirsthead
\multicolumn{5}{r}{\textit{(The table continues in the next page)}}
\endfoot
\multicolumn{5}{l}{Table 3  continue}
\endhead

\endlastfoot

\hline
\rule[-1 mm]{0mm}{5 mm}
    Size  & Type   &
  \begin{minipage}{2.3 cm}
 \smallskip
  \centering
  List of \\
 points \\
 \medskip
  \end{minipage}   &
  \begin{minipage}{2.3 cm}
 \smallskip
  \centering
  Stabilizer\\
 group \\
 \medskip
  \end{minipage}   &
  \begin{minipage}{0.9 cm}
 \smallskip
  \centering
 Cons-\\
truc-\\
 tion \\
 \medskip
  \end{minipage}   \\
   \hline
\rule[-1 mm]{0mm}{5 mm}
$ 5 $ & CC & $ \{1, 2, 4, 8, 15\}  $ & $ S_5 $ &  B\\
\hline
\rule[-1 mm]{0mm}{5 mm}
$ 6 $ &  NC & $ \{1, 2, 3, 4, 8, 12\}  $ & $ (S_3 \times S_3) \rtimes C_2 $ & A  \\
\hline
\rule[-1 mm]{0mm}{5 mm}
$ 8 $ & CC & $ \{1, 2, 4, 7, 8, 11, 13, 14\}  $ & $  ASL(3,2) $ &  H \\
\hline
\rule[-1 mm]{0mm}{5 mm}
$ 8 $ & NC & $ \{1, 2, 4, 5, 8, 9, 12, 13\}  $ & $ PSL(3,2) $ &  A \\
\hline
  \end{longtable}
   \end{footnotesize}

\smallskip

\newpage


   \begin{footnotesize}
   \begin{longtable}
 {|c|c|c|c|c|}
 \caption{Classification of the minimal 1-saturating sets in $\PG(4,2)$}

\endfirsthead
\multicolumn{5}{r}{\textit{(The table continues in the next page)}}
\endfoot
\multicolumn{5}{l}{ \textbf{Table 5}  continue}
\endhead

\endlastfoot
   \hline
\rule[-1 mm]{0mm}{5 mm}
    Size  & Type   &
  \begin{minipage}{2.3 cm}
 \smallskip
  \centering
  List of \\
 points \\
 \medskip
  \end{minipage}   &
  \begin{minipage}{2.3 cm}
 \smallskip
  \centering
  Stabilizer\\
 group \\
 \medskip
  \end{minipage}   &
    \begin{minipage}{0.9 cm}
 \smallskip
  \centering
 Cons-\\
truc-\\
 tion \\
 \medskip
  \end{minipage}   \\
   \hline
\rule[-1 mm]{0mm}{5 mm}

$ 9 $ & CC & $ \{1, 2, 4, 8, 14, 16, 22, 27, 28\}  $ & $ C_2 \times PSL(3,2) $ &  $L_{21}$ \\
\hline
\rule[-4 mm]{0mm}{10 mm}
$ 9 $ & NC & $ \{1, 2, 4, 6, 8, 16, 20, 22, 27\}  $ & $ S_3 \times S_4 $ & 
\begin{minipage}{0.4 cm}
B,\\
GL
\end{minipage}
\\
\hline
\rule[-1 mm]{0mm}{5 mm}

$ 10 $ & CC & $ \{1, 2, 4, 8, 15, 16, 21, 22, 27, 28\}  $ & $  G_{1920} $   & D  \\
\hline
\rule[-1 mm]{0mm}{5 mm}

$ 10 $ & NC & $ \{1, 2, 4, 5, 8, 10, 16, 22, 27, 28\}  $ & $ D_8 $ &   \\
\hline
\rule[-1 mm]{0mm}{5 mm}
$ 10 $ & NC & $ \{1, 2, 4, 8, 10, 16, 20, 22, 23, 27\}  $ & $ D_{12} $ &   \\
\hline
\rule[-1 mm]{0mm}{5 mm}
$ 10 $ & NC & $ \{1, 2, 4, 8, 10, 14, 16, 17, 22, 28\}  $ & $ C_2 \times S_4 $ &   \\
\hline
\rule[-1 mm]{0mm}{5 mm}
$ 10  $ & NC & $ \{1, 2, 4, 5, 8, 11, 16, 22, 27, 28\}  $ & $ C_2 \times S_4 $ &   \\
\hline

  \end{longtable}
   \end{footnotesize}


   \begin{footnotesize}
   \begin{longtable}
 {|c|c|c|c|c|}

  \caption{Classification of the minimal 1-saturating sets in $\PG(5,2)$}
\endfirsthead
\multicolumn{5}{r}{\textit{(The table continues in the next page)}}
\endfoot
\multicolumn{5}{l}{ Table 5  continue}
\endhead

\endlastfoot

\hline
    Size  & Type   &
  \begin{minipage}{2.3 cm}
 \smallskip
  \centering
  List of \\
 points \\
 \medskip
  \end{minipage}   &
  \begin{minipage}{2 cm}
 \smallskip
  \centering
  Stabilizer\\
 group \\
 \medskip
  \end{minipage}   &
 \begin{minipage}{0.9 cm}
 \smallskip
  \centering
 Cons-\\
truc-\\
 tion \\
 \medskip
  \end{minipage}   \\
\hline
\rule[-2 mm]{0mm}{8 mm}   	
$ 13 $ & CC & $  \begin{minipage}{5 cm}   \{1, 2, 4, 7, 8, 16, 25, 32, 37, 38, 43, 51, 58\}   \end{minipage}   $ & $  G_{1152} $ & $L_{21}$  \\
\hline
\rule[-4 mm]{0mm}{10 mm}
$ 13 $ & NC & $  \begin{minipage}{5 cm}   \{1, 2, 4, 7, 8, 14, 16, 20, 25, 32, 43, 52, 63\}   \end{minipage}   $ &
\begin{minipage}{2 cm}
$ (C_2 \times C_2 \times C_2 \times C_2) \rtimes C_2 $
\end{minipage}
&   \\
\hline
\rule[-2 mm]{0mm}{8 mm}
$ 13 $ & NC & $  \begin{minipage}{5 cm}   \{1, 2, 4, 7, 8, 16, 20, 24, 25, 32, 37, 43, 46\}   \end{minipage}   $ & $ C_2 \times S_4 $ &   \\
\hline
\rule[-2 mm]{0mm}{8 mm}
$ 13 $ & NC & $  \begin{minipage}{5 cm}   \{1, 2, 4, 7, 8, 16, 20, 25, 29, 32, 37, 43, 46\}   \end{minipage}   $ & $ C_2 \times C_2 \times S_4 $ &   \\
\hline
\rule[-2 mm]{0mm}{8 mm}
$ 13 $ & NC & $  \begin{minipage}{5 cm}   \{1, 2, 4, 7, 8, 16, 17, 24, 25, 32, 37, 38, 43\}   \end{minipage}   $ & $ C_2 \times C_2 \times S_4 $ & GL \\
\hline
\rule[-2 mm]{0mm}{8 mm}
$ 13 $ & NC & $  \begin{minipage}{5 cm}   \{1, 2, 4, 7, 8, 9, 16, 21, 23, 32, 33, 58, 59\}   \end{minipage}   $ & $ C_2 \times C_2 \times S_4 $ &   \\
\hline
\rule[-2 mm]{0mm}{8 mm}
$ 13 $ & NC & $  \begin{minipage}{5 cm}   \{1, 2, 4, 5, 6, 7, 8, 16, 27, 32, 43, 48, 59\}   \end{minipage}   $ & $  G_{1152} $ &  B \\
\hline

  \end{longtable}
   \end{footnotesize}

$\;$\\
$\;$\\

   \begin{footnotesize}

   \begin{longtable}
 {|c|c|c|c|c|}

 \caption{Classification of small minimal 1-saturating sets in $\PG(6,2)$}

\endfirsthead
\multicolumn{5}{r}{\textit{(The table continues in the next page)}}
\endfoot
\multicolumn{5}{l}{ Table 6  continue}
\endhead

\endlastfoot

   \hline
\rule[-2 mm]{0mm}{8 mm}
    Size  & Type   &
  \begin{minipage}{1.9 cm}
 \smallskip
  \centering
  List of \\
 points \\
 \medskip
  \end{minipage}   &
  \begin{minipage}{2.1 cm}
 \smallskip
  \centering
  Stabilizer\\
 group \\
 \medskip
  \end{minipage}   &
  \begin{minipage}{0.9 cm}
 \smallskip
  \centering
 Cons-\\
truc-\\
 tion \\
 \medskip
  \end{minipage}   \\
   \hline
\rule[-4 mm]{0mm}{10 mm}
$ 19 $ & NC &
$  \begin{minipage}{5.8 cm}
 \{1, 2, 4, 8, 15, 16, 31, 32, 43, 51, 55, \\
 64, 67, 85, 89, 101, 110, 121, 126\}
  \end{minipage}
 $ &
\begin{minipage}{2 cm}
 $ (C_2 \times C_2 \times C_2 \times C_2) \rtimes C_2 $
\end{minipage}
&   \\
\hline
\rule[-2 mm]{0mm}{8 mm}
$ 19 $ & NC & $  \begin{minipage}{5.8 cm}   \{1, 2, 4, 8, 15, 16, 26, 29, 32, 39, \\ 43, 51, 64, 70, 76, 85, 110, 120, 121\}   \end{minipage}  $ & $ S_5 $ &   \\
\hline
\rule[-4 mm]{0mm}{10 mm}
$ 19 $ & NC & $  \begin{minipage}{5.8 cm}   \{1, 2, 4, 6, 8, 14, 15, 16, 24, 32,  \\ 43, 47, 48, 50, 51, 64, 85, 108, 121\}   \end{minipage}  $ &
\begin{minipage}{2 cm}
\begin{center}
$ (C_2 \times C_2 \times A_5) \rtimes C_2 $
\end{center}
\end{minipage}
&   \\
\hline
\rule[-2 mm]{0mm}{8 mm}
$ 19 $ & NC & $  \begin{minipage}{5.8 cm}   \{1, 2, 4, 8, 13, 15, 16, 22, 30, 32,  \\ 42, 43, 48, 51, 55, 64, 85, 108, 121\}   \end{minipage}  $ &
\begin{minipage}{2 cm}
\begin{center}
$ (A_4 \times A_5) \rtimes C_2 $
\end{center}
\end{minipage}
&   \\
\hline
\rule[-2 mm]{0mm}{8 mm}
$ 19 $ & NC & $  \begin{minipage}{5.8 cm}   \{1, 2, 4, 8, 15, 16, 30, 32, 43, 51,  \\ 54, 64, 66, 85, 89, 101, 108, 120, 127\}   \end{minipage}  $ & $  G_{5760} $ &   \\
\hline
\rule[-2 mm]{0mm}{8 mm}
$ 20 $ & NC & $  \begin{minipage}{5.8 cm}   \{1, 2, 4, 8, 15, 16, 29, 31, 32, 37,  43,  \\ 51, 64, 72, 85, 99, 110, 118, 121, 126\}   \end{minipage}  $ & $ C_2 \times C_2 $ &   \\
\hline

  \end{longtable}
   \end{footnotesize}

\section*{Acknowledgements}
The research of S. Marcugini and F. Pambianco was supported in
part by Ministry for Education, University and Research of Italy (MIUR) (Project "Geometrie
di Galois e strutture di incidenza") and by the Italian National Group for Algebraic
and Geometric Structures and their Applications (GNSAGA - INDAM).
 The research of A.A. Davydov was carried out at the IITP RAS at the expense of the
Russian Foundation for Sciences (project 14-50-00150).

\end{document}